\documentclass[mathematics,article,accept,pdftex,moreauthors]{Definitions/mdpi} 
\firstpage{1} 
\pubvolume{1}
\issuenum{1}
\articlenumber{0}
\pubyear{2024}
\copyrightyear{2024}
\externaleditor{\textls[-25]{Academic Editor: Mario Versaci}}
\datereceived{24 July 2024} 
\daterevised{10 August 2024} 
\dateaccepted{13 August 2024} 
\datepublished{ } 
\hreflink{https://doi.org/} 
\usepackage{color}

\definecolor{c1}{rgb}{0,0,1} 

\Title{A Mathematical Programming Model for Minimizing Energy Consumption on a Selective Laser Melting {Machine}}

\TitleCitation{A Mathematical Programming Model for Minimizing Energy Consumption on a Selective Laser Melting Machine}

\Author{Chunlong 
 Yu * and Junjie Lin}

\AuthorNames{Chunlong Yu and Junjie Lin}

\AuthorCitation{Yu, C.; Lin, J.}

\address[1]{%
School of Mechanical Engineering, Tongji University, Shanghai 201804, China.} 

\corres{\hangafter=1 \hangindent=1.05em \hspace{-0.82em} Correspondence: chunlong\_yu@tongji.edu.cn}

\abstract{The scheduling problem in additive manufacturing is receiving increasing attention; however, few have considered the effect of scheduling decisions on machine energy consumption. This research focuses on the nesting and scheduling problem of a single selective laser melting (SLM) machine to reduce total energy consumption. Based on an energy consumption model, a nesting and scheduling problem is formulated, and a mixed integer linear programming model is proposed. This model simultaneously determines part-to-batch assignments, part placement in the batch, and the choice of build orientation to reduce the total energy consumption of the SLM machine. The energy-saving potential of the model is validated through numerical experiments. Additionally, the effect of the number of alternative build orientations on energy consumption is explored.}



\begin{document}
	
	
	
\section{Introduction}
Compared with conventional manufacturing, additive manufacturing (AM) has great potential to improve material efficiency, reduce lifecycle impact, and~enable greater engineering capabilities~\cite{huang2016energy}. Despite these advantages, AM is recognized as an energy-intensive technology, especially for large-scale parts and mass production. The~formation of multiple thin layers from raw materials requires a considerable amount of energy, resulting in higher energy consumption per unit volume of material compared to conventional manufacturing techniques \citep{yoon2014comparison}. This aspect raises concerns about the potential negative economic and environmental impacts. Therefore, reducing energy consumption is becoming one of the most important goals driving the adoption of additive manufacturing in key industries \citep{karimi2021energy}.

The research of energy-aware scheduling in conventional manufacturing systems is receiving considerable attention. Various studies have considered the multi-objective scheduling problem and proposed methods to optimize the energy cost with other common scheduling criteria~\cite{jia2017bi, zhou2018multi, wang2018bi}. In~these studies, the~energy-saving opportunities come from a better allocation of jobs to machines with different speeds and powers, as~well as a better job sequence to avoid high consumption in the periods of high electricity~price. 

However, few works have considered the scheduling problem that aims to reduce energy consumption in AM machines. In~this research, we focus on saving the energy consumption of the selective laser melting (SLM) machine. SLM is a powder-bed-fusion AM process whereby a high-density-focused laser beam selectively scans a powder bed and those scanned and solidified layers are stacked upon each other to build a fully functional three-dimensional part, tool, or prototype~\cite{kruth2005benchmarking}. SLM machines are batch processing machines that are capable of printing a group of parts at a time. According to a recent study by Lv~et~al.~\cite{lv2021novel}, the~energy consumption of an SLM machine is related to the power of each subsystem and the subsystem running time. The~power of each subsystem depends on the machine parameters, which are basically constant during the process. The~running time of each subprocess depends on not only the process parameters (laser power, heating and cooling temperature, etc.), but~also the batch parameters (total part volume, surface area, support structure volume, height, etc.). Therefore, the~energy consumption can be affected by how the batches are formed. This creates an opportunity to optimize the energy of printing the parts by making better decisions such as part-to-batch assignment, nesting, and build orientation~selection.

Generally, increasing the batch utilization reduces the total energy consumption. Compared to a single part, printing a full batch at a time can save more than 50\% energy consumption per part~\cite{faludi2017environmental}. To~maximize the batch utilization, one needs to consider the nesting problem, which decides the part displacement and rotation around the z axis (vertical axis) and packs parts into batches. For~SLM machines, no stacks are permitted because supporting structures are required for heat dissipation and fixing parts on the platform~\cite{zhang2020improved}. This results in a 2D nesting problem that concerns packing the projections of the parts into the build platform of the machine, which is usually a rectangle.  To~reduce the problem difficulty, a~common way is to represent the geometry of a part with a bin. However, even for this case, the~nesting problem is NP-hard~\cite{lodi2002recent}.

\begin{figure}[H]
	\includegraphics[width=10 cm]{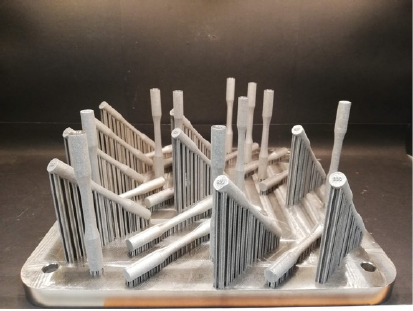}
	\caption{A batch of parts with different build~orientations.\label{fig:Batch_diff_ori}}
\end{figure}

On the other hand, a~part generally has several optional build orientations that can satisfy the quality requirement. Changing the build orientation can lead to a variation in the part height, support structure volumes, and~the projection area of the part, as~shown in Figure~\ref{fig:Batch_diff_ori}. The~batch height determines the number of layers and thus the total recoating time of a batch, and~the volume of the support structure determines the scanning time of the support. These two factors have an impact on the energy consumption of several subsystems of the SLM machine. Therefore, the~choice of the build orientation is also important for energy saving. At~the same time, a~proper selection of the build orientation can increase the packing efficiency and thus reduce the energy~consumption.

In this study, we consider the nesting and scheduling problem that jointly handles the part-to-batch assignment, part nesting and build orientation selection on a single SLM machine to minimize the total energy consumption. The~contributions of this research~are as follows: 
\begin{itemize}
	\item	The energy-saving potential of nesting and scheduling decisions on a single SLM machine is studied. Based on an accurate energy consumption model, we propose a mixed-integer linear programming (MILP) model considering alternative build orientations for the nesting and scheduling problem;
	\item	Through numerical experiments, the~performance of the proposed model is tested and validated by comparison with commercial 3D printing software. The~MILP model achieves significant energy savings compared to the commercial software. The~percentage of energy savings tends to increase with the number of alternative build orientations of the parts.
\end{itemize}

What remains of this paper is organized as follows. Section \ref{sec2} reviews the literature. \mbox{Section \ref{sec3}} describes the problem. Section \ref{sec4} describes the energy consumption model. \mbox{Section~\ref{sec5}} proposes the mathematical model. Section \ref{sec6} reports the numerical results, and \mbox{Section \ref{sec7}} concludes the~paper.

\section{Literature~Review}\label{sec2}

In this section, we first review the works related to the energy-aware scheduling, then we introduce some works about nesting and scheduling on AM machines, and~finally we summarize the research~gap.

\subsection{Energy-Aware~Scheduling}
Many scholars have paid attention to the problem of energy-aware scheduling in manufacturing systems with conventional machines. Agrawal~et~al.~\cite{agrawal2014energy} proposed three~energy-aware scheduling algorithms including genetic algorithm, cellular automata and efficiency-based allocation heuristic to schedule tasks on a non-identical machines environment to minimize the makespan and energy consumption. Bruzzone~et~al.~\cite{bruzzone2012energy} integrated the energy-aware scheduling module with the advanced pand scheduling system to control the power peaks on the shop floor given a detailed schedule. Che~et~al.~\cite{che2017energy} studied the scheduling problem of unrelated parallel machines with the goal of minimizing the total power cost under a scenario that electricity prices may vary from hour to hour throughout a day.  Schulz~et~al.~\cite{schulz2019multi} used a local search algorithm to solve a multi-objective problem of energy consumption, energy cost, and~demand charging in a flow shop. Che~et~al.~\cite{che2016efficient} developed a continuous-time MILP model and proposed an efficient greedy insertion heuristic to minimize the total electricity cost within a given makespan under time-of-use or time-dependent electricity tariffs. In~addition, Abikarram~et~al.~\cite{abikarram2019energy} addressed a parallel machine scheduling problem, in~which the number of active machines and utilization are controlled to reduce the peak~power.

For additive manufacturing, there are only a few related studies on energy-aware scheduling. Most studies aim to improve the energy efficiency of AM machines by optimizing the process parameters to fabricate parts with low energy consumption while meeting desired specifications. For~instance, Paul~et~al.~\cite{paul2012process} considered the problem of finding the orientation of the part and the layer thickness that yield the minimum energy consumption for selective laser sintering (SLS). Verma~et~al.~\cite{verma2017sustainability} proposed a multi-objective optimization model to find the combination of processing parameters that minimizes laser processing energy and material waste while maximizing part quality for a single SLS process. For~the first time, Karimi~et~al.~\cite{karimi2021energy} proposed an MILP model that integrates process-level controls and system-level scheduling decisions. By~adjusting the process parameters at different printing stages and controlling the starting time of each job, the~model successfully reduces peak power and, thus, the~cost of electricity tariffs for on-demand~charging.

\subsection{Nesting and Scheduling on AM~Machines}
Scheduling problems for AM machines can be categorized into three types: nesting, scheduling, nesting and scheduling~\cite{oh2020nesting, pinto2024nesting}. In~the past, nesting and scheduling were usually discussed separately. Recently, the~number of articles focusing on the joint problem of nesting and scheduling is increasing. Chergui~et~al.~\cite{chergui2018production} addressed the nesting and scheduling problem with the due date constraints. They proposed a heuristic approach to minimize the maximum lateness. Alicastro~et~al.~\cite{alicastro2021reinforcement} considered a nesting and scheduling problem to minimize the makespan of unrelated parallel machines. The~authors proposed a local search algorithm enhanced by a reinforcement learning mechanism for neighborhood selection. In~this research, the~build orientations of the parts are prefixed. Li~et~al.~\cite{li2018single} studied a scheduling problem in a batch processing machine with 2D capacity constraint to minimize the makespan. They proposed a greedy heuristic and a genetic algorithm to solve the problem. Che~et~al.~\cite{che2021machine} studied an unrelated parallel AM machine scheduling problem that requires simultaneously assigning parts to batches, orienting parts, packing parts onto a 2D platform, and~assigning batches to machines to minimize the makespan. An~efficient simulated annealing approach was proposed. Yu~et~al.~\cite{yu2022mathematical} investigated the nesting and scheduling problem in an unrealted parallel machines environment to minimize the total tardiness of orders. They proposed different mathematical models and concluded that simplifying the complex nesting features of the problem can improve the solvability of the problem, but~also lead to the potential packing infeasibility of batches. Zipfel~et~al.~\cite{zipfel2024iterated} considered the scenario where sequence-dependent setups are required between batches of different material. They proposed an iterated local search metaheuristic to minimize the total weighted tardiness of orders. Recently, Nascimento~\cite{nascimento2024optimal} proposed an exact solution approach based on Bender's decomposition framework and constrain programming to solve the nesting and scheduling of parts with complex~shapes. 

\subsection{Research~Gap}
Many studies have focused on reducing the energy consumption and costs by scheduling or adjusting the machine on--off state. But~few have paid attention to the energy-saving benefit brought by nesting and scheduling decisions in AM~machines.

Baumers et al. \cite{baumers2011energy} compared the energy consumption of a build with a single part to that of multiple parts, and the~full-build utilization produced a large (97.79\%) energy saving on the EOSINT P 390 printer. They also showed that printing a single part in different build orientations consumes different amounts of energy. Furthermore, Piili~et~al.~\cite{piili2015cost} found that a full build can reach a cost reduction of 81--92\% compared to single part build, and~the nesting solution is the major variable that users can influence the energy in a multi-part fabrication environment. Both studies reveal the potential of energy saving in AM machines by a proper nesting and scheduling, but~a solution to obtain the optimal nesting and scheduling decisions is not~provided. 

On the other hand, although~the nesting and scheduling problem has been investigated in several studies, none have considered the objective of energy savings except for our conference paper~\cite{lin2023bi}, which emphasizes the trade-off between productivity and energy consumption in multiple~machines. 

{In summary, our work is driven by the exploration of energy-saving potential in SLM machines through optimized nesting and scheduling decisions. The~primary objective is to formulate the nesting and scheduling problem as a mathematical programming model. By~solving this model, one can determine the optimal part-to-batch assignments, part nesting configurations, and~build orientation selections when printing a set of parts on a single SLM machine.}

\section{Problem~Description}\label{sec3}
Consider a set of parts \( J = \{1, \ldots, n\} \) to be printed on a single SLM machine. Each part \( j \in J \) has a given volume \( v_j \) and surface \( s_j \) and a set of optional build orientations \( K_j \). Each build orientation \( k \in K_j \) corresponds to a 3D geometry with four parameters \( \{h_{jk}, l_{jk}, w_{jk}, s_{jk}\} \), which are the height, length, and width of the rectangle projection, and~the volume of the supporting structure, respectively. The~SLM machine can process a batch of parts at a time. The~capacity of the machine is represented by a cuboid of length \( LW \), width \( WW \), and height \( HW \). A~group of parts can be placed on the same batch if their projections do not overlap and their geometries are within the platform boundaries. Parts must be placed parallel to the length or width of the platform and can be rotated by 90° around the Z axis. The~goal is to decide the build orientation for each part, place and schedule them on the machines to minimize the total energy~consumption.

{We impose the following~assumptions:
	\begin{enumerate}
		\item For any part $j \in J$, any optional build orientation $k \in K_j$ can satisfy the quality requirement of the part $j$.
		
		\item Each part is represented by its minimum bounding box, and~thus, the~projection of any part is a rectangle on the build~platform.
		
		\item There is a minimum distance requirement between any parts placed on the build platform. Without~loss of generality, it is assumed that this distance has been included in the size of each~part.
		
		\item The process parameters (e.g., laser power, layer thickness, hatching distance) are identical for all~parts.
		
		\item The power of each subsystem of the SLM machine is assumed to be constant during a~subprocess.
	\end{enumerate}
}

\section{Energy Consumption~Model}\label{sec4}
Energy forecasting is the key for the evaluation and reduction in SLM energy consumption. Energy consumption models for AM can be divided into three categories: specific energy models~\cite{kellens2017environmental}, stage-based energy models~\cite{faludi2017environmental}, and subsystem-based energy models~\cite{yi2018energy}. 

Models that compute the energy demand by multiplying the specific energy consumption (SEC) and the weight or volume of a part are called specific energy models. Here, the~SEC refers to the energy consumed for each deposited kilogram or volume of material. It may be affected by the utilization of the build volume, build density, and orientation, and hence, different values were reported in the literature for the same AM technology~\cite{ingarao2018environmental, yi2018energy}. To~describe the energy consumption in a more reasonable way, substage-based energy consumption models were proposed. These models predict the total energy consumption by summing the energy consumed in each stage. The~power and time of each operation stage are obtained through measurement, and~thus, these models can only be applied on the specific AM machine. Subsystem-based models work in a similar way, they estimate the total energy by summing the energy consumed by each machine subsystem. In~this research stream, most works focus on the energy modeling of the machine subsystem~\cite{paul2012process}, while the operating status of different subsystems at different stages are~ignored.

Recently, Lv~et~al.~\cite{lv2021novel} proposed a more accurate energy consumption model that jointly considers the power of the machine subsystem, the~time of the subprocess, and~the working state of the machine subsystem in each subprocess. In~this paper, we adopt this model and combine it with the   mathematical programming model for energy~optimization.

\subsection{Energy Consumption~Calculation}

\begin{figure}[H]
	\includegraphics[width=12 cm]{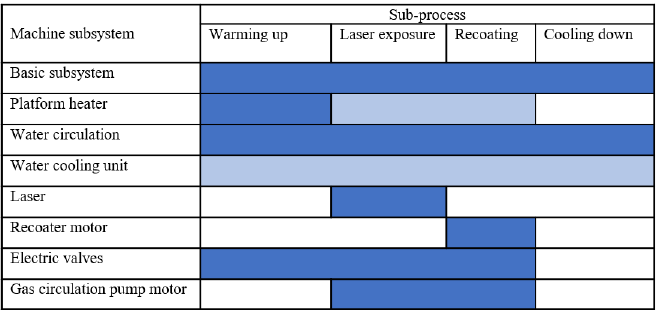}
	\caption{Working status of different machine subsystems during the SLM process. The~dark color, light color, and blank areas represent the subsystems running all the time, intermittently, and out of~service. \label{fig:SLM_states}}
\end{figure}

Typically, an~SLM machine consists of eight subsystems, and~the process of printing a batch can be divided into four subprocesses, as~shown in Figure~\ref{fig:SLM_states}. The~energy consumption of the SLM machine for printing a single batch is calculated as:
\begin{equation}
	E = \sum_{f=1}^{F} E_f = \sum_{f=1}^{F} \sum_{l=1}^{L} p_f t_l k_{fl} = PKT,
	\label{eq:energy_equation}
\end{equation}
where \( E_f \) is the energy consumption of the \( f \)-th subsystem, \( p_f \) is the power of the \( f \)-th subsystem, \( t_l \) is the time of the \( l \)-th subprocess, and~\( k_{fl} \in [0,1] \) is a factor describing the working state of the \( f \)-th subsystem in the \( l \)-th subprocess. In~the following subsections, we describe the way to calculate the processing time of subprocesses, the~power of subsystems, and the state of each subsystem in each~subprocess.

\subsection{Processing Time of~Subprocesses}
A typical SLM process includes four subprocesses: heating, laser exposure, recoating, and~cooling. During~the heating phase, the~following relationship holds:
\begin{equation}
	t_h = C_1 T^2 + C_2 T - C_3, \label{eq:t_heating_exp}
\end{equation}
where $t_h$ represents the time spent in the heating stage, $T$ represents the temperature of the platform at a certain moment. The~parameters $C_1$, $C_2$, and $C_3$ are machine-dependent factors and are obtained by regression using the experimental data. See Lv~et~al.~\cite{lv2021novel} for more technical details. Thus, the~time of the heating phase is expressed as:
\begin{equation}
	\Delta t_h = t_h(T_f) - t_h(T_i),
\end{equation}
\textls[-25]{where $T_i$ is the initial temperature of the platform, and~$T_f$ is the final temperature \mbox{after~heating.}}

The construction stage is the main production process, including powder recoating and laser exposure. The~total time spent in this phase is expressed as:
\begin{equation}
	t_b = t_l + t_r,
\end{equation}
where $t_l$ is the laser exposure time and $t_r$ is the powder spreading time. The~time of the sub-stages of the laser exposure process is expressed as:
\begin{equation}
	t_l = t_b + t_c + t_p + t_s = \frac{\sum_{j \in B} s_j}{n_l v_b \Delta y} + \frac{\sum_{j \in B} s_j}{n_l v_c \Delta y} + \frac{\sum_{j \in B} V_j}{v_p} + \frac{\sum_{j \in B} s_{j k}}{v_s},
\end{equation}
where $t_b$, $t_c$, $t_p$, $t_s$ are the laser exposure time to scan the border, fill contour, hatch volume, and~build supports, $s_j$ is the area of the outer and inner surfaces of part $j$ scanned as boundaries during construction [mm\(^2\)], $B$ is the set of parts in the batch, $n_l$ is the number of lasers operating during manufacturing, $v_b$ and $v_c$ are the speeds for scanning the boundaries and filling contours [mm/s], respectively; $\Delta y$ is the layer thickness [mm], $V_j$ and $s_{jk}$ are the volumes of the part and support [mm\(^3\)], respectively; $\dot{v_p}$ and $\dot{v_s}$ are the build-up rate for the part and supports [mm\(^3\)/s], both can be expressed as:
\begin{equation}
	\dot{v} = n_l \times D \times \Delta y \times v,
\end{equation}
where $D$ is the hatching distance [mm/s], $v$ is the scanning speed for part or support. The~recoating time can be expressed as:
\begin{equation}
	t_r = N t_{r0},
\end{equation}
where $t_{r0}$ is the time required for spreading a layer of powder. The~number of slices depends on the height of the batch and the layer thickness, which can be expressed as:
\begin{equation}
	N = \frac{H}{\Delta y}.
\end{equation}

After the last layer has been scanned, the~platform is gradually cooled down in the inert gas used to prevent the oxidation of the metal powder. The~following relationship holds:
\begin{equation}
	t_c = C_4 T^2 - C_5 T + C_6,
\end{equation}
where $T$ is the temperature of the platform at a certain moment, and the~constants $C_4$, $C_5$, and $C_6$ are obtained by regression from the experimental data. Hence, the~time of the cooling phase is expressed as:
\begin{equation}
	\Delta t_c = t_c(T_c) - t_c(T_b), \label{eq:t_cooling}
\end{equation}
where $T_b$ is the temperature of the platform before cooling and $T_c$ is the final~temperature.

\subsection{Power of~Subsystems}

The power parameters of each subsystem are summarized in Table~\ref{table:power_of_subsystems}. The~values of each parameter are machine-dependent and can be measured by the power acquisition device, including data acquisition chassis, data acquisition card, voltage sensor, current sensor, etc. For~details, please refer to Lv~et~al.~\cite{lv2021novel}. Note that the power of the laser subsystem could be different during the scanning of the border, filling contour, hatching volume, and building supports; for~convenience, we treat the laser used in these subprocesses as four different subsystems, and thus the total energy consumption can be obtained by Equation~(\ref{eq:energy_equation}). 

\begin{table}[H] 
	\caption{Power of the different~subsystems.\label{table:power_of_subsystems}}
	\newcolumntype{C}{>{\centering\arraybackslash}X}
	\begin{tabularx}{\textwidth}{LL}
		\toprule
		\textbf{Subsystem} & \textbf{Power [W]}\\
		\midrule
		Basic subsystem & $P_{bs}$\\
		Platform heater & $P_{ht}$\\
		Water circulation unit & $P_{wc}$\\
		Water-cooling unit & $P_{co}$\\
		Laser-scanning border & $P_{lsb}$\\
		Laser-filling contour & $P_{lfc}$\\
		Laser volume hatching & $P_{lvh}$\\
		Laser support structure & $P_{lss}$\\
		Recoater motor & $P_{rc}$\\
		Electric valves & $P_{ev}$\\
		Gas circulation pump motor & $P_{gp}$\\
		\bottomrule
	\end{tabularx}
\end{table}

\subsection{State of~Subsystems}
Generally, a~subsystem of an SLM machine can be fully on, off, or~running intermittently in different subprocesses. For~example, the~platform heater is fully on during the preheating phase to raise the temperature of the platform to the required level, while during the laser exposure phase, it runs intermittently to keep the temperature constant. Therefore, to~correctly estimate the energy consumption of each subsystem, it is critical to identify their actual operating states during each~subprocess.

To this end, the~working state factor $k_{fl}$ is introduced. For~the $f$-th subsystem, $k_{fl}$ is the ratio of the actual running time to the total time of the $l$-th subprocess. As~shown in Figure~\ref{fig:SLM_states}, the~basic subsystem and water circulation unit are running all the time, whilst the laser, powder spreading motor, electric valve, and~gas circulation pump motor only operate at certain stages. In~these cases, the~$k_{fl}$ is either 1 or 0. The~platform heater and water cooling unit operate intermittently. To~decide the values of $k_{fl}$ for these two subsystems, measurements are often required. See Lv~et~al.~\cite{lv2021novel} for details of the measurement~procedure.
\section{Mathematical Programming~Model}\label{sec5}
In this section, we propose an MILP model for the defined problem and several means to improve the model efficiency. After, we introduce some improvements upon the~model.

\subsection{MILP~Model}
The MILP model is given as follows. The~sets~are  as follows:
\begin{itemize}
	\item $J$: set of~parts.

	\item $F = \{ \textit{bs}, \textit{ht}, \textit{wc}, \textit{co}, \textit{lsb}, \textit{lfc}, \textit{lvh}, \textit{lss}, \textit{rm}, \textit{ev}, \textit{gp} \}$: set of subsystems. \textit{bs} is the basic system, \textit{ht} is the heater, \textit{wc} is the water circulation system, \textit{co} is the water cooling unit, \textit{lsb} is the laser unit for scanning border, \textit{lfc} is the laser unit for filling contour, \textit{lvh} is the laser unit for volume hatching, \textit{lss} is the laser unit for scanning supporting structure, \textit{rm} is the recoater motor, \textit{ev} stands for the electric values, \textit{gp} is the gas circulation pump~motor.
	
	\item $L = \{ \textit{ph}, \textit{sb}, \textit{fc}, \textit{vh}, \textit{ss}, \textit{rc}, \textit{co} \}$: set of subprocesses. \textit{ph} is the preheating, \textit{sb} is the scanning border, \textit{fc} is the filling contour, \textit{vh} is volume hatching, \textit{ss} is the building of support structure, \textit{rc} is the recoating, \textit{co} is the~cooling.
	
	\item $K_j$: set of alternative build orientations of part $j$.
	\item $B$: set of available batches on machine.
\end{itemize}

The parameters~are  as follows:
\begin{itemize}
	\item $v_j$: volume of part $j$.
	\item $a_j$: surface area of part $j$.
	\item $D$: the minimum allowed distance from the parts to the platform boundary.
	\item $d$: the minimum allowed distance between parts.
	\item $t_{r0}$: recoating speed of machine.
	\item $LW, WW, HW$: length, width, and~height of the build platform.
	\item $h_{jk}, l_{jk}, w_{jk}, s_{jk}$: height, length, width, and~supporting structure volume of the $k$-th orientation of part $j$.
	\item $n_l$: number of lasers.
	\item $v_b$: speed for scanning border.
	\item $v_c$: speed for filling contour.
	\item $\dot{v_p}$: build-up rate for the part.
	\item $\dot{v_s}$: build-up rate for the support.
	\item $\Delta y$: layer thickness.
	\item $p_f$: power of the subsystem $f$.
	\item $TPH$: the time spent in the preheating subprocess.
	\item $TCO$: the time spent in the cooling subprocess.
	\item $\phi_{fl}$: the power coefficient of the subsystem $f$ in the subprocess $l$, $\phi_{fl} = 0$ means that the subsystem $f$ is not active in the subprocess $l$, $\phi_{fl} = 1$ means that the subsystem $f$ is fully active during the subprocess $l$; $\phi_{fl} \in (0,1)$ means that the subsystem $f$ is intermittently running during the subprocess $l$, that is, $\phi_{fl} = \frac{P_{fl}^*}{P_f}$ where $P_{fl}^*$ represents the measured power of subsystem $f$ in subprocess $l$.
\end{itemize}

The decision variables~are as follows:
\begin{itemize}
	\item $X_{jb}$: equals 1 if part $j$ is assigned to the $b$-th batch, 0 otherwise;
	\item $Y_{jb}$: equals 1 if part $j$ selects the $k$-th optional build orientation, 0 otherwise;
	\item $(x_j, y_j)$: coordinates of the left-bottom point of part $j$'s projection on the platform;
	\item $O_j$: equals 1 if part $j$ is placed such that its length is parallel to that of the platform.
\end{itemize}

The auxiliary variables~are as follows:
\begin{itemize}
	\item $Z_b$: equals 1 if the $b$-th batch on machine is formed, 0 otherwise;
	\item $PL_{jj'}$: equals 1 if part $j$'s right--top point is placed to the left of part $j$'s left--bottom point in the same batch, 0 otherwise;
	\item $PB_{jj'}$: equals 1 if part $j$'s right--top point is placed below part $j$'s left--bottom point in the same batch, 0 otherwise;
	\item $H_b$: height of the $b$-th batch;
	\item $E_b$: the energy consumption of the $b$-th batch;
	\item $T_{lb}$: the processing time of the subprocess $l$ in the $b$-th batch;
	\item $M$: a large number.
\end{itemize}

In the model, the~decision variables $PL_{jj'}$, $PB_{jj'}$, $PL_{j'j}$, $PB_{j'j}$ together describe the positional relationship of parts $j$ and $j'$. 

The nesting and scheduling problem of a single SLM machine is formulated as:
\begin{equation}
	\min \sum_{b \in B} E_b, \label{eq:obj}
\end{equation}
\begin{equation}
	\sum_{j \in J} X_{jb} \leq M Z_b, \quad \forall b \in B, \label{eq:c1}
\end{equation}
\begin{equation}
	\sum_{b \in B} X_{jb} = 1, \quad \forall j \in J,  \label{eq:c2}
\end{equation}
\begin{equation}
	\sum_{k \in K_j} Y_{jk} = 1, \quad \forall j \in J,  \label{eq:c3}
\end{equation}
\begin{equation}
	\sum_{j \in J} \frac{a_j X_{jb}}{n_l v_b \Delta y} \leq T_{sb,b}, \quad \forall b \in B,  \label{eq:c4}
\end{equation}
\begin{equation}
	\sum_{j \in J} \frac{a_j X_{jb}}{n_l v_c \Delta y} \leq T_{fc,b}, \quad \forall b \in B,  \label{eq:c5}
\end{equation}
\begin{equation}
	\sum_{j \in J} \frac{v_j X_{jb}}{\dot{v_p}} \leq T_{vh,b}, \quad \forall b \in B,  \label{eq:c6}
\end{equation}
\begin{equation}
	\sum_{j \in J} \frac{X_{jb} \sum_{k \in K_j} s_{jk} Y_{jk}}{\dot{v_s}} \leq T_{ss,b}, \quad \forall b \in B,  \label{eq:c7}
\end{equation}
\begin{equation}
	T_{ph,b} = TPH, \quad \forall b \in B,  \label{eq:c8}
\end{equation}
\begin{equation}
	T_{co,b} = TCO, \quad \forall b \in B,  \label{eq:c9}
\end{equation}
\begin{equation}
	H_b \geq \sum_{k \in K_j} Y_{jk} h_{jk} - M (1 - X_{jb}), \quad \forall b \in B, j \in J,  \label{eq:c10}
\end{equation}
\begin{equation}
	t_{r0}\frac{H_b}{\Delta y} \leq T_{rc,b}, \quad \forall b \in B,  \label{eq:c11}
\end{equation}
\begin{equation}
	E_b \geq \sum_{f \in F} \sum_{l \in L} p_f \phi_{fl} T_{lb}, \quad \forall b \in B,  \label{eq:c12}
\end{equation}
\begin{equation}
	H_b \leq HW, \quad \forall b \in B,  \label{eq:c13}
\end{equation}
\begin{equation}
	x_j + \sum_{k \in K_j} Y_{jk} l_{jk} \leq LW + M (1 - O_j), \quad \forall j \in J,  \label{eq:c14}
\end{equation}
\begin{equation}
	x_j + \sum_{k \in K_j} Y_{jk} w_{jk} \leq LW + MO_j, \quad \forall j \in J,  \label{eq:c15}
\end{equation}
\begin{equation}
	y_j + \sum_{k \in K_j} Y_{jk} w_{jk} \leq WW + M (1 - O_j), \quad \forall j \in J,  \label{eq:c16}
\end{equation}
\begin{equation}
	y_j + \sum_{k \in K_j} Y_{jk} l_{jk} \leq WW + MO_j, \quad \forall j \in J,  \label{eq:c17}
\end{equation}
\begin{equation}
	x_j + \sum_{k \in K_j} Y_{jk} l_{jk} + d \leq x_{j'} + M (1 - PL_{jj'}) + M (1 - O_j), \quad \forall j, j' \in J,  \label{eq:c18}
\end{equation}
\begin{equation} 
	x_j + \sum_{k \in K_j} Y_{jk} w_{jk} + d \leq x_{j'} + M (1 - PL_{jj'}) + MO_j, \quad \forall j, j' \in J, \label{eq:c19}
\end{equation}
\begin{equation}
	y_j + \sum_{k \in K_j} Y_{jk} w_{jk} + d \leq x_{j'} + M (1 - PB_{jj'}) + M (1 - O_j), \quad \forall j, j' \in J, \label{eq:c20}
\end{equation}
\begin{equation}
	y_j + \sum_{k \in K_j} Y_{jk} l_{jk} + d \leq x_{j'} + M (1 - PB_{jj'}) + MO_j, \quad \forall j, j' \in J, \label{eq:c21}
\end{equation}
\begin{equation}
	PL_{jj'} + PB_{jj'} + PL_{j'j} + PB_{j'j} \geq X_{jb} + X_{j'b} - 1, \quad \forall j, j' \in J, j < j', b \in B, \label{eq:c22}
\end{equation}
\begin{equation}
	Z_b \in \{0,1\}, \quad \forall b \in B, \label{eq:c23}
\end{equation}
\begin{equation}
	T_{lb} \geq 0, \quad \forall l \in L, b \in B, \label{eq:c24}
\end{equation}
\begin{equation}
	X_{jb} \in \{0,1\}, \quad \forall j \in J, b \in B, \label{eq:c25}
\end{equation}
\begin{equation}
	Y_{jk} \in \{0,1\}, \quad \forall j \in J, k \in K_j, \label{eq:c26}
\end{equation}
\begin{equation}
	H_b \geq 0, \quad \forall b \in B, \label{eq:c27}
\end{equation}
\begin{equation}
	PL_{jj'}, PB_{jj'} \in \{0,1\}, \quad \forall j, j' \in J, \label{eq:c28}
\end{equation}
\begin{equation}
	x_j, y_j \geq 0, \quad \forall j \in J, \label{eq:c29}
\end{equation}
\begin{equation}
	O_j \in \{0,1\}, \quad \forall j \in J, \label{eq:c30}
\end{equation}

\textls[-31]{Constraints (\ref{eq:c1}) ensure that the batches cannot be assigned without forming. \mbox{Constraints (\ref{eq:c2})}} specify to assign parts to only one batch. Constraints (\ref{eq:c3}) guarantee that the part selects only one build orientation. Constraints (\ref{eq:c4}) calculate the time for the laser to scan the border. Constraints (\ref{eq:c5}) calculate the time for the laser to fill contour. Constraints (\ref{eq:c6}) calculate the time for the laser to hatch the volume. Constraints (\ref{eq:c7}) calculate the time for the laser to build the support structure. Constraints (\ref{eq:c8}) define the preheating time. Constraints (\ref{eq:c9}) define the cooling time. Constraints (\ref{eq:c10}) define the batch height. Constraints (\ref{eq:c11}) calculate the powder recoating time. Constraints (\ref{eq:c12}) calculate the energy consumption of the batch. Constraints~(\ref{eq:c13}) guarantee that the height of the batch is less than the height of the machine. Constraints~(\ref{eq:c14})--(\ref{eq:c17}) ensure that parts do not exceed the machine boundary. Constraints~(\ref{eq:c18})--(\ref{eq:c21}) ensure that when parts are assigned to the same batch, they will not overlap. Specifically, when $PL_{jj'}$ equals 1, part $j$ is placed on the left side of part $j'$, so the right edge of part $j$ does not exceed the left edge of part $j'$. When $PL_{jj'}$ equals 0, this constraint is released. Constraints (\ref{eq:c20}) and (\ref{eq:c21}) express the same relationship in the vertical direction. Constraints (\ref{eq:c22}) guarantee that when part $j$ and part $j'$ are assigned to the same batch, i.e.,~$X_{jb}$, $X_{j'b}$ equal to 1, part $j$ and $j'$ choose at least one of the positional relationships. For~example, when $PL_{jj'}$ equals 1, $PB_{jj'}$ equals 0, $PL_{j'j}$ equals 0 and $PB_{j'j}$ equals 0, part $j'$ is on the left semi-axis of part $j$. When $PL_{jj'}$ equals 1, $PB_{jj'}$ equals 1, $PL_{j'j}$ equals 0, and $PB_{j'j}$ equals 0, part $j'$ is located in the third quadrant of part $j$. Another positional relationship can be obtained by setting different values of $PL_{jj'}$ and $PB_{jj'}$. 

\subsection{Model~Improvements}
We improve the model as follows. First, the~following constraints are introduced to break the symmetry of batches:
\begin{equation}
	Z_{(b-1)} \leq Z_b, \quad \forall b \in B. \label{eq:c31}
\end{equation}

Second, 
 to~linearize (\ref{eq:c7}), we introduce an auxiliary variable $e_{jb}$ to represent the support structure volume of part $j$ in the $b$-th batch, and~replace Constraints  (\ref{eq:c7}) with the following:
\begin{equation}
	e_{jb} \geq \sum_{k \in K_j} Y_{jk} s_{jk} - M(1 - X_{jb}), \quad \forall j \in J, b \in B, \label{eq:c32}
\end{equation}
\begin{equation}
	\sum_{j \in J} \frac{e_{jb}}{v_s} \leq T_{ss,b}, \quad \forall b \in B. \label{eq:c33}
\end{equation}

Finally, in~order to further improve the solution efficiency of the model, we define the upper bound of the big-$M$ values in different constraints, as~reported in Table~\ref{table:Big-M}. 

\begin{table}[H]
	\caption{Big-M values for different~constraints.}
	\centering
	\begin{tabularx}{\textwidth}{LL}
		\toprule
		\textbf{Constraints} & \textbf{Big-M value} \\
		\midrule
		(\ref{eq:c1}) & $M = |J|$ \\
		(\ref{eq:c10}) & $M = \max_{j \in J, k \in K_j} h_{jk}$ \\
		(\ref{eq:c14}) (\ref{eq:c17}) & $M = \max_{j \in J, k \in K_j} l_{jk}$ \\
		(\ref{eq:c15}) (\ref{eq:c16}) & $M = \max_{j \in J, k \in K_j} w_{jk}$ \\
		(\ref{eq:c18}) & $M = \max_{j \in J, k \in K_j} l_{jk} + LW$ \\
		(\ref{eq:c19}) & $M = \max_{j \in J, k \in K_j} w_{jk} + LW$ \\
		(\ref{eq:c20}) & $M = \max_{j \in J, k \in K_j} w_{jk} + WW$ \\
		(\ref{eq:c21}) & $M = \max_{j \in J, k \in K_j} l_{jk} + WW$ \\
		(\ref{eq:c32}) & $M = \max_{j \in J, k \in K_j} s_{jk}$\\
		\bottomrule
	\end{tabularx}
	\label{table:Big-M}
\end{table}

\subsection{Model Parameter: The Number of Available~Batches}
{Besides the big-$M$ values, another critical parameter influencing model efficiency is the number of available batches on the machine, i.e.,~how many batches can be opened on the machine. Let $n_B$ denote the cardinality of $B$. A~large $n_B$ results in a substantial number of decision variables $X_{jib}$, which may complicate the resolution process. Conversely, a~small $n_B$ may lead to suboptimal solutions when $n_B$ is smaller than the number of batches required in the optimal solution. In~the worst case, the~model could become infeasible if $n_B$ is set below the minimum number of batches necessary to accommodate all parts. Thus, determining an appropriate value for $n_B$ is~crucial.
	
	Let $n_B^{\text{fea}}$ represent the minimum number of batches required to ensure model feasibility, and~let $n_B^{\text{opt}}$ denote the number of batches in the optimal solution. It is evident that
\begin{equation}
		n_B^{\text{fea}} \leq n_B^{\text{opt}} \leq |J|,
	\end{equation}
	where $|J|$ is the number of parts to be printed. However, $n_B^{\text{fea}}$ and $n_B^{\text{opt}}$ are unknown and difficult to determine precisely. To~address this, we propose a simple trial-and-error approach as~follows:
\begin{enumerate}[label=\text{Step \arabic*:}]
	\item Set 
 $n_B = \left\lceil \eta |J| \right\rceil$, where $\eta \in (0,1)$;
	\item Solve the MILP model with the current value of $n_B$;
	\item If the model is infeasible or if the obtained solution includes $n_B$ opened batches, then increment $n_B$ by 1 and return to Step 2. Otherwise, proceed to the next step;
	\item Terminate and output the results.
\end{enumerate}
}

\section{Numerical~Results}\label{sec6}
\subsection{Test~Instance}
The part models are downloaded from the website \url{https://www.stlfinder.com/3dmodels/free-3d-models/} (accessed on 13 August 2024). All the part information such as volume, surface, support structure, and~height are obtained from the STL files using a commercial 3D printing software named Materialise Magics. {Details of the parts are provided in Appendix \ref{Appendx:A}. }

\subsection{Machine and Process~Parameters}
The parts are to be fabricated using aluminum AlSi10Mg powder on an SLM280HL machine. The~machine is equipped with two 400 W fiber lasers with a laser beam diameter of 80 µm. {The length, width, and height of the build platform of the machine is 268 mm, 268 mm, and 315 mm, respectively.  The~power of machine subsystems and the power coefficients of the subsystems in different subprocesses are summarized in Table~\ref{Table: Power_and_coefficient}. }

The parts are to be fabricated in an inert gas atmosphere. Before~the construction process began, the~oxygen content will be reduced from 21\% to below 0.1\% to prevent the oxidation and reaction of the AlSi10Mg powder. The~platform will be heated from 27~°C to 150 °C to relieve thermally induced internal stress of the fabricated part, which results in a preheating time of $\Delta t_c = t_c (80~^\circ \text{C}) - t_c (150~^\circ \text{C})$. During~the manufacturing process, two lasers work together. The~time required for spreading a layer of powder is 11~s. After~building the parts, the~platform will be cooled from 150~°C to 80~°C, which leads to a cooling time of $\Delta t_c = t_c (80~^\circ \text{C}) - t_c (150~^\circ \text{C})$. The~process parameters of the SLM machine are given in Table~\ref{table:ProcessParams}.

\begin{table}[H] 
	\caption{{The} power and power coefficient of the different~subsystems. \label{Table: Power_and_coefficient}}
	\newcolumntype{L}{>{\raggedright\arraybackslash}m{3.45cm}}
	\newcolumntype{R}{>{\raggedleft\arraybackslash}m{1.35cm}}
	\newcolumntype{C}{>{\raggedleft\arraybackslash}m{1.45cm}}
	
	\tablesize{\footnotesize}
\begin{adjustwidth}{-\extralength}{0cm}

\begin{tabularx}{\fulllength}{LRCCCCCCc}
		\toprule
		\textbf{Subsystem} & \textbf{Power [W]} & \textbf{Pre Heating} & \textbf{Scan Border} & \textbf{Fill Contour} & \textbf{Volume Hatching} & \textbf{Support Structure} & \textbf{Re Coating} & \textbf{Cooling} \\
		\midrule
		Basic subsystem & $569.7$ & $1$ & $1$ & $1$ & $1$ & $1$ & $1$ & $1$\\
		Platform heater & $1122.3$ & $1$ & $0.4826$ & $0.4826$ & $0.4826$ & $0.4826$ & $0.4826$ & $0$\\
		Water circulation unit & $713.3$ & $1$ & $1$ & $1$ & $1$ & $1$ & $1$ & $1$\\
		Water-cooling unit & $1739.4$ & $0.168$ & $0.353$ & $0.353$ & $0.353$ & $0.353$ & $0.353$ & $0.216$\\
		Laser-scanning border & $1770.9$ & $0$ & $1$ & $0$ & $0$ & $0$ & $0$ & $0$\\
		Laser-filling contour & $1770.9$ & $0$ & $0$ & $1$ & $0$ & $0$ & $0$ & $0$\\
		Laser volume hatching & $2022.9$ & $0$ & $0$ & $0$ & $1$ & $0$ & $0$ & $0$\\
		Laser support structure & $2022.9$ & $0$ & $0$ & $0$ & $0$ & $1$ & $0$ & $0$\\
		Recoater motor & $52.1$ & $0$ & $0$ & $0$ & $0$ & $0$ & $1$ & $0$\\
		Electric valves & $32.1$ & $1$ & $1$ & $1$ & $1$ & $1$ & $1$ & $0$\\
		Gas circulation pump motor & $69.1$ & $0$ & $1$ & $1$ & $1$ & $1$ & $1$ & $0$\\
		\bottomrule
	\end{tabularx}
\end{adjustwidth}
\end{table}
\unskip

\begin{table}[H]
	\caption{Process~parameters.}
	\centering
	\begin{tabularx}{\textwidth}{lCCcc}
		\toprule
		\textbf{Process Parameter} & \textbf{Scan Border} & \textbf{Fill Contour} & \textbf{Volume Hatching} & \textbf{Support Structure}\\
		\midrule
		\text{Power of the laser [W]} 
 & 300 & 300 & 350 & 350\\
		\text{Scanning speed [mm/s]} & 730 & 730 & 1650 & 1000\\
		\text{Layer thickness [mm]} & 0.03 & 0.03 & 0.03 & 0.03\\
		\text{Hatching distance [mm]} & - & - & 0.13 & 0.13\\
		\bottomrule
	\end{tabularx}
	\label{table:ProcessParams}
\end{table}
\unskip

\subsection{Experiment~Setup}
With the above settings, we conducted two~experiments:
\begin{enumerate}
	\item \text{Experiment 1}: 
 The first numerical experiment compares the energy consumption of printing a set of 20 parts when MILP and Magics are applied for the nesting and scheduling. Magics uses the default build orientation for each part, which corresponds to the smallest Z axis height. For~MILP, we assume that each part has five different alternative build orientations, which are provided by the Magics software according to different rules. The~MILP model is solved by the Gurobi solver (ver 9.5.2) with a maximum CPU runtime of~7200 s. 
	
	\item \text{Experiment 2}: A design of experiments is performed with the following factors and levels: number of parts---\{20, 25, 30\}, number of alternative build orientations---\{1, 3, 5, 7\}. This results in 12 experiments. We analyze the impact of the number of parts and alternative build orientations on the energy consumption of the MILP model.
\end{enumerate}
\subsection{Results~Analysis}
\subsubsection{Experiment~1}

\begin{figure}[H]
	\includegraphics[width=11 cm]{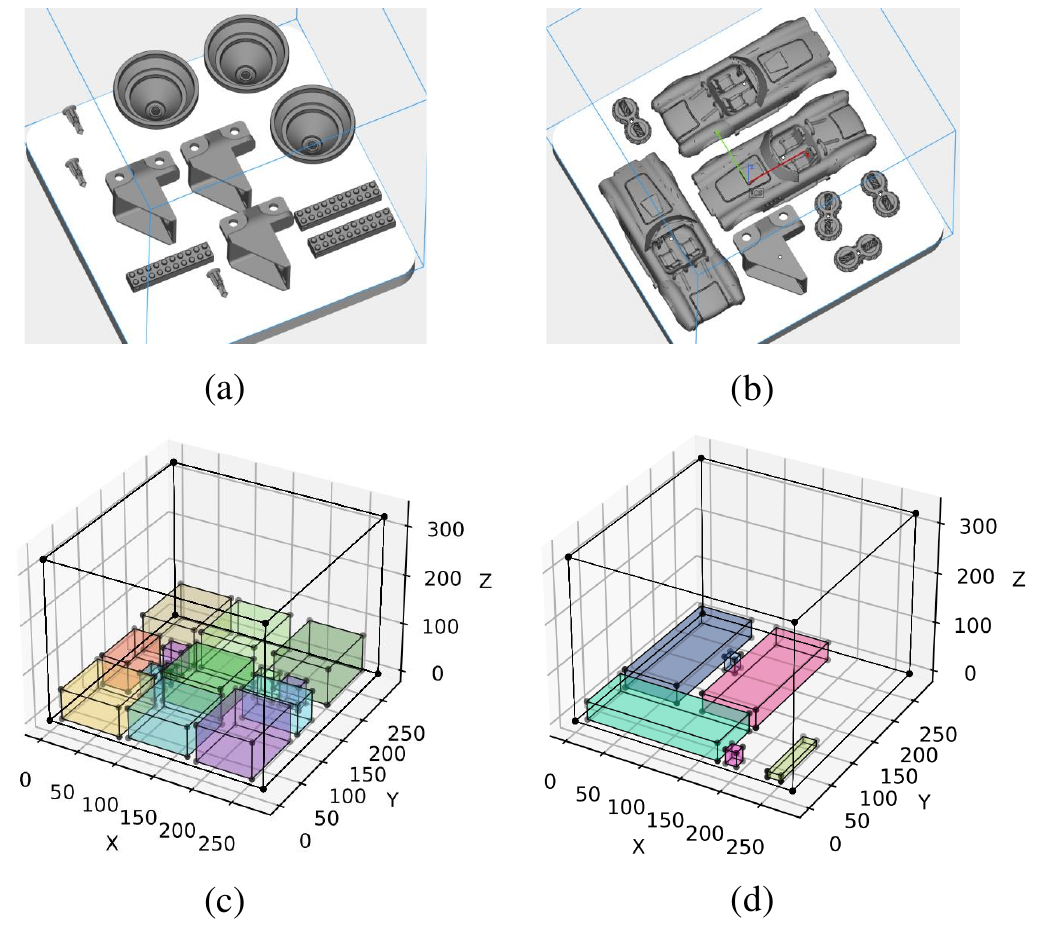}
	\caption{The placement for the scheduling by MILP and Magics on ec\_20-5: (\textbf{a}) Magics Batch 1; (\textbf{b})~Magics Batch 2; (\textbf{c}) MILP Batch 1; and (\textbf{d}) MILP Batch~2. \label{fig:CompPlacement}}
\end{figure}

The batch placements of the Magics and MILP are illustrated in Figure~\ref{fig:CompPlacement}. In~addition, we also summarize the batch information, printing time, and~energy consumption of the two approaches in Table~\ref{table:Comparison_MILP_Magics}. Compared to Magics, the~MILP model reduces the energy consumption by 41.67 MJ (8.62\%). Furthermore, the~total processing time of MILP is less than that of Magics. Indeed, the~energy consumption of the machine is determined by both the time of each subprocess and the power of each subsystem. In~the case of constant energy power, the~reduction in energy consumption can only be achieved by reducing the running time of~subprocesses. 

\vspace{-9pt}

\begin{table}[H]
	\caption{Comparison of the results of two~methods.}
	\centering
\tablesize{\small}	
\begin{adjustwidth}{-\extralength}{0cm}

\begin{tabularx}{\fulllength}{lCCC CCC}
		\toprule
		\text{Parameters} & \multicolumn{3}{c}{\textbf{Magics}} & \multicolumn{3}{c}{\textbf{MILP}} \\
		\cmidrule(lr){2-4} \cmidrule(lr){5-7}
		& \text{Batch 1} & \textbf{Batch 2} & \textbf{Total} & \textbf{Batch 1} & \text{Batch 2} & \textbf{Total} \\
		\midrule
		\text{Surface area [$\text{mm}^2$]} 
 & 220,689 & 382,340 & 603,029 & 276,929 & 326,100 & 603,029 \\
		\text{Part volume [$\text{mm}^3$]} & 146,217 & 188,334 & 334,551 & 186,356 & 148,195 & 334,551 \\
		\text{Support volume [$\text{mm}^3$]} & 108,969 & 141,951 & 250,920 & 43,618 & 112,625 & 156,243 \\
		\text{Number of slices} & 2031 & 1729 & 3760 & 2481 & 1220 & 3701 \\
		\text{Utilization rate} & 47.63\% & 57.84\% & 52.7\% & 63.72\% & 47.70\% & 55.7\% \\
		\text{Number of parts} & 12 & 8 & 20 & 14 & 6 & 20 \\
		\text{Makespan [s]} & 63,749 & 77,965 & 141,714 & 68,851 & 63,448 & 132,299 \\
		\text{Total energy consumption [MJ]} & 225.66 & 296.59 & 522.25 & 241.84 & 238.72 & 480.56 \\
		\bottomrule
	\end{tabularx}
\end{adjustwidth}
	\label{table:Comparison_MILP_Magics}
\end{table}

Table~\ref{table:energy_saving_distribution} presents the energy savings associated with different subsystems. The~``Percentage'' column indicates the proportion of energy savings contributed by each subsystem, calculated as follows:
\begin{equation}
\frac{E_f (\text{MILP}) - E_f(\text{Magic})} {\text{Total energy saving}} \times 100, 
\end{equation}
where $E_f (\text{MILP})$ is the total energy consumed by subsystem $f$ in the MILP solution, and $E_f (\text{Magics})$ is that of the Magics solution. This table shows that the laser for the support structure provides the greatest energy savings, followed by the water circulation unit, water cooling unit, basic subsystem, and~platform heater. In~contrast, the~recoater motor, electric valves, and~gas circulation pump contribute relatively little to energy savings. This is because the MILP solution, which involves less subprocess time, leads to more energy savings in subsystems with a higher power consumption. Thus, reductions in the subprocess time have a greater impact on energy-intensive units, explaining why the support structure construction shows significant energy savings when there is a notable difference in the time required for building the support structure in the solutions of MILP and Magics. {Similarly, Table~\ref{table:energy_saving_distribution_process} presents the energy savings associated with different subprocesses. As~shown, the~majority of energy savings is achieved during the support structure building subprocess, with~the remainder obtained during the recoating phase.}

\begin{table}[H]
	\caption{The energy saving contribution of different~subsystems.  }
	\centering
	\begin{tabularx}{\textwidth}{lCC}
		\toprule
		\textbf{Machine Subsystem}
 & \textbf{Energy Saving [MJ]} & \textbf{Percentage [\%]}\\
		\midrule
		\text{1-Basic subsystem}  & 5.36 & 12.86\\
		\text{2-Platform heater} & 5.1 & 12.24\\
		\text{3-Water circulation unit} & 6.72 & 16.13\\
		\text{4-Water cooling unit} & 5.78 & 13.87\\
		\text{5-Laser for support structure} & 17.73 & 42.55\\
		\text{6-Recoater motor} & 0.03 & 0.07\\
		\text{7-Electric valves} & 0.3 & 0.72\\
		\text{8-Gas circulation pump motor} & 0.65 & 1.56\\
		\midrule
		\text{Total energy saving} & 41.67 & 100\\
		\bottomrule
	\end{tabularx}
	\label{table:energy_saving_distribution}
\end{table}
\unskip

\begin{table}[H]
	\caption{The energy saving contribution of different~subprocess. }
	\centering
	\begin{tabularx}{\textwidth}{LCC}
		\toprule
		\textbf{Machine Subprocess} & \textbf{Energy Saving [MJ]} & \text{Percentage [\%]}\\
		\midrule
		\text{1-Preheating} & 0 & 0\\
		\text{2-Scan border} & 0 & 0\\
		\text{3-Fill contour} & 0 & 0\\
		\text{4-Volumne hatching} & 0 & 0\\
		\text{5-Support structure} & 40.00 & 95.99\\
		\text{6-Recoating} & 1.67 & 4.01\\
		\text{7-Cooling} & 0 & 0\\
		\midrule
		\text{Total energy saving} & 41.67 & 100\\
		\bottomrule
	\end{tabularx}
	\label{table:energy_saving_distribution_process}
\end{table}

\begin{figure}[H]
	\includegraphics[width=11 cm]{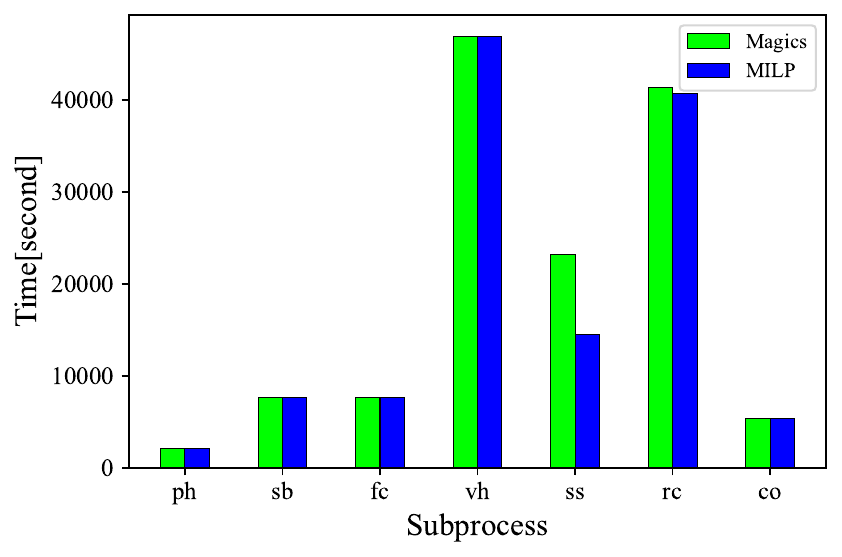}
	\caption{The duration of different subprocesses. \textit{Ph} is preheating, \textit{sb} is scan border, \textit{vh} is volume hatching, \textit{ss} is support structure, \textit{rc} is powder spreading, and \textit{co} is~cooling.  \label{fig:CompTime}}
\end{figure}

\begin{figure}[H]
	\includegraphics[width=11 cm]{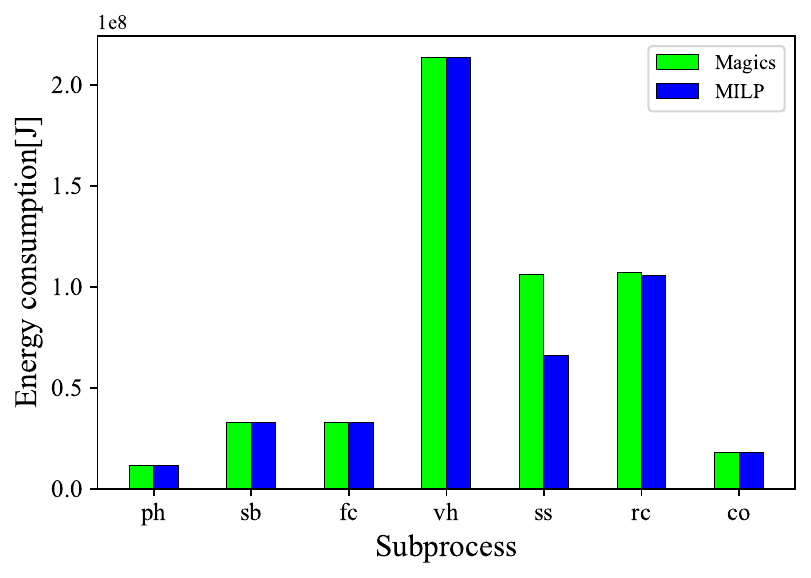}
	\caption{The energy consumption of different subprocesses.  \textit{Ph}
 is preheating, \textit{sb} is scan border, \textit{vh} is volume hatching, \textit{ss} is support structure, \textit{rc} is powder spreading, and \textit{co} is cooling. 
		\label{fig:CompEnergeProcess}}
\end{figure}

\begin{figure}[H]
	\includegraphics[width=11 cm]{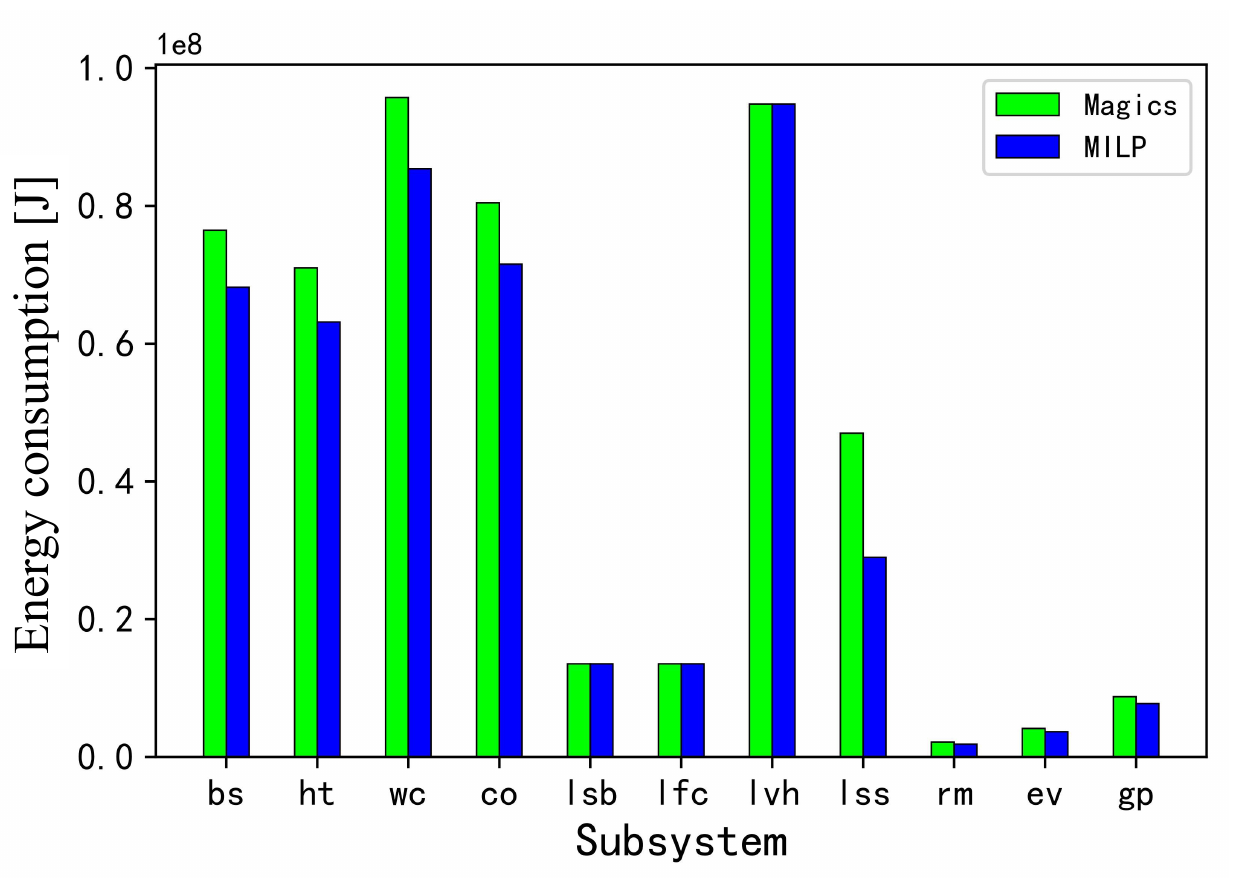}
	\caption{The energy consumption of different subsystems. \textit{bs} is basic subsystem, \textit{ht} is platform heater, \textit{wc} is water circulation unit, \textit{co} is water cooling unit, \textit{lsb} is laser for scanning border, \textit{lfc} is laser for filling contour, \textit{lvh} is laser for volume hatching, \textit{lss} is laser for support structure, \textit{rm} is Recoater, \textit{ev} is electric valves, and \textit{gp} is gas circulation pump motor.
		\label{fig:CompEnergy}}
\end{figure}

To investigate the reasons behind energy savings and time reductions, we compare the two solutions obtained by Magics and MILP in terms of time differences across different subprocesses (Figure \ref{fig:CompTime}), {energy consumption across different subprocesses (Figure \ref{fig:CompEnergeProcess})}, and~energy consumption for different subsystems (Figure \ref{fig:CompEnergy}). As~illustrated in Figure~\ref{fig:CompTime}, the~primary differences between the two solutions are observed in the time required to build the support structure and the time needed to spread the powder. According to formulas (\ref{eq:t_heating_exp})--(\ref{eq:t_cooling}), for~the same printer and process parameters, heating and cooling times are fixed. Additionally, when the part volume and surface area are identical, the~times for scanning the border, filling contour, and~volume hatching remain constant. The~time required to print the support structure is dependent on its volume, while the powder recoating time is influenced by the number of slices. Table~\ref{table:Comparison_MILP_Magics} shows that variations in the support structure volume and the number of slices between the two solutions account for the differences in processing time. Furthermore, Figure~\ref{fig:SLM_states} reveals that nearly all subsystems are active during both the laser exposure and recoating stages. Since the laser system and recoater motor operate alternately during these stages, changes in the time required to build the support structure and recoating time significantly impact the energy consumption of most subsystems. This explains the observed reductions in energy~consumption.

Overall, the~primary factors contributing to energy savings are the reductions in time during the recoating phase and the construction of the support structure. Consequently, the printing efficiency improves while optimizing energy consumption. For~a nesting problem, the~choice of the part build orientation significantly impacts both the processing time and energy consumption of the batch. This is because the build orientation affects the support structure, the~height of the part, and~the projection area of the part, which in turn influences batch~utilization.

\unskip
\unskip

\subsubsection{Experiment~2}
The results of instances with different a number of parts and build orientations are summarized in Table~\ref{table:instance_comparison}, where \text{ins\_20\_1}  represents the instance having 20 parts and each part has only one build orientation 1, \text{ins\_20\_3} represents the instance with 20 parts and each has three different build orientations, etc. The~column ``EC'' represents the total energy consumption, the~column ``EV'' represents the energy consumption per unit volume, and~the column ``ES'' represents the proportion of energy-saving of MILP compared to Magics. Specifically, the~energy consumption of Magics when printing 20 parts, 25 parts, and~30 parts is 522.25 MJ, 663.79 MJ, and~883.27 MJ, respectively. The~column ``Gap'' is the optimality gap of the best feasible solution. The~Column ``B'' represents the number of batches opened in the solution. All the EC values are obtained from the MILP~model.

\begin{table}[H]
	\caption{Energy consumption with different numbers of parts and build~orientations.}
	\centering
	\begin{tabularx}{\textwidth}{CCCCCC}
		\toprule
		\textbf{Instance} & \textbf{EC [MJ]} & \textbf{EV [J/mm$^3$]} & \textbf{ES [\%]} & \textbf{Gap [\%]} & \textbf{B}\\
		\midrule
		\text{ins\_20\_1} & 507.73 & 841.97 & 3.46 & 0.00 & 2\\
		\text{ins\_20\_3} & 481.06 & 797.74 & 8.53 & 5.17 & 2\\
		\text{ins\_20\_5} & 480.56 & 796.91 & 8.62 & 5.34 & 2\\
		\text{ins\_20\_7} & 479.91 & 795.83 & 8.75 & 6.46 & 2\\
		\text{ins\_25\_1} & 625.23 & 839.38 & 5.81 & 7.66 & 3\\
		\text{ins\_25\_3} & 602.57 & 808.96 & 9.22 & 5.32 & 3\\
		\text{ins\_25\_5} & 577.11 & 774.78 & 13.06 & 1.58 & 2\\
		\text{ins\_25\_7} & 570.31 & 765.65 & 14.08 & 1.59& 2\\
		\text{ins\_30\_1} & 782.87 & 840.81 & 11.37 & 7.40 & 3\\
		\text{ins\_30\_3} & 746.03 & 801.25 & 15.54 & 12.29 & 3\\
		\text{ins\_30\_5} & 745.24 & 800.40& 15.63 & 15.76 & 3\\
		\text{ins\_30\_7} & 769.41 & 826.36 & 12.89 & 20.21& 3\\
		\bottomrule
	\end{tabularx}
	\label{table:instance_comparison}
\end{table}

\begin{figure}[H]
	\includegraphics[width=11 cm]{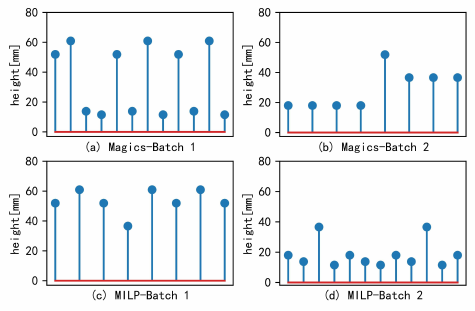}
	\caption{Height distribution of parts across different batches for two scheduling solutions from instance ins\_20\_1.
		\label{fig:PartHeightDistribution}}
\end{figure}

It can be observed from Table~\ref{table:instance_comparison} that even when all parts share a single build orientation, consistent with Magics' build direction, MILP can still achieve some energy savings. In~this scenario, the~volume, surface area, and~support structure volume of all parts remain identical in both solutions; the only variation lies in part-to-batch allocation, specifically in how the parts are grouped into batches. Figure~\ref{fig:PartHeightDistribution} illustrates the differences in part height within the batches provided by Magics and MILP in instance ins\_20\_1. In~Batch 1 of both solutions, the~batch height (maximum part height within the batch) is the same at 60.9 mm. However, in~Batch 2, the~batch height in Magics is 51.9 mm, whereas MILP results in a lower batch height of 36.6 mm. This indicates that, even with the same number of parts and build orientations, different scheduling results can lead to varying total batch heights and, consequently, different energy~consumption.

Additionally, when the number of printed parts remains the same, increasing the build orientations of the parts can also reduce the energy consumption. For~example, with~20 parts, as~the number of build orientations increases from 1 to 7, the~energy consumption per unit volume of parts gradually decreases, and~the energy savings achieved by MILP are also greater. It is understandable that increasing the number of build orientations for a part can reduce energy consumption. Different build orientations result in varying support structure volumes and part heights. By~choosing a build orientation with a smaller support structure volume without increasing the batch height, energy consumption can be significantly reduced. Furthermore, as~the number of parts increases, the~advantage of the MILP model over Magics becomes more pronounced. For~instance, when the number of printed parts increases from 20 to 25, as~seen in instances \text{ins\_20\_5} 
 and \text{ins\_25\_5}, the~percentage of energy savings rises significantly. This trend is also observed in other instances, where increased numbers of printed parts with the same build orientation result in greater energy savings compared to Magics~scheduling.

To summarize, the~following observations can be obtained from Table~\ref{table:instance_comparison}:
\begin{enumerate}
	\item Increasing the build orientation of parts can reduce energy consumption to some extent, but~it also complicates the model and increases the optimality gap.
	\item Increasing the number of printed parts can lead to reduced energy consumption compared to Magics, but~it also makes the model more challenging to solve.
\end{enumerate}

\vspace{-6pt}
\unskip

\section{Conclusions}\label{sec7}
In this study, we address the nesting and scheduling problem on a single SLM machine with the goal of reducing our energy consumption. We consider the 2D nesting subproblem involved in packing parts into batches and account for the multiple optional build orientations for each part. We propose a mathematical programming model and several enhancements to strengthen~it.

By comparing the results from solving the MILP model with those provided by commercial 3D printing software, we demonstrate the potential for energy savings through optimized nesting and scheduling. {The energy savings range from 3.46\% to 15.63\% across various instances, depending on the number of parts and build orientations.} Specifically, when the build orientation of parts is fixed, the~MILP model reduces the total height of batches by optimizing the distribution of parts among batches. This reduces the total number of layers and, consequently, the~overall energy consumption. When multiple build orientations are available, the~MILP model achieves even greater energy savings due to the added flexibility in optimization. However, this increased flexibility also raises the complexity of solving the model, necessitating a careful~trade-off.

The numerical results indicate that there is still potential for further energy savings through improvements in model efficiency. Future work will focus on developing efficient solution approaches, such as matheuristics and metaheuristics, to~handle large-scale problems. {Additionally, this study addresses the nesting and scheduling problem on a single SLM machine. However, since multi-machine environments better reflect industrial practices, extending the research to such environments will also be considered. This extension would require accounting for various factors, including different additive manufacturing technologies, build platform sizes, and~the powder materials used by the machines on the shop floor.
}

\authorcontributions{Conceptualization, C.Y.; methodology, C.Y.; software, J.L.; validation, C.Y. and J.L.; formal analysis, C.Y.; investigation, J.L.; resources, C.Y.; data curation, J.L.; writing---original draft preparation, J.L.; writing---review and editing, C.Y.; visualization, J.L.; supervision, C.Y.; project administration, C.Y.; funding acquisition, C.Y. All authors have read and agreed to the published version of the~manuscript.}

\funding{This research was funded by National Natural Science Foundation of China (72301196), and~Tongji University ``Fundamental Research Funds for the Central Universities''.}



\dataavailability{The raw data supporting the conclusions of this article and the source codes will be made available by the authors upon request.} 

\conflictsofinterest{The authors declare no conflicts of interest.
	
}

\appendixtitles{yes} 
\appendixstart
\appendix
\section[\appendixname~\thesection]{Part Information}
\label{Appendx:A}
{The information of each part type is given in Table~\ref{table:PartInfo}. Each part can have maximum optional build orientations which are provided by the Materialise Magics software according to different rules. We used in total six different part types to generate the instances. For~example, ins\_20\_5 contains four copies of part 1 and part 2, and~three copies of part 3 to part 6. Each part has five optional build orientations 1--5. }

\begin{table}[H]
	\caption{Part~parameters. \label{table:PartInfo}}
	\tablesize{\footnotesize}
	\begin{adjustwidth}{-\extralength}{0cm}
		\newcolumntype{C}{>{\centering\arraybackslash}X}
		\begin{tabularx}{\fulllength}{ccccCccc}
			\toprule
			\textbf{Part Type} & \textbf{Build Orientation} & \textbf{Volume [mm$^3$]} & \textbf{Surface Area [mm$^2]$} & \textbf{Supporting Structure Volume~[mm$^3$]} & \textbf{Length [mm]} & \textbf{Width [mm]} & \textbf{Heigth  [mm]} \\
			\midrule
			\text{1}
 & 1 & 6744.00 & 8607.80 & 1724.00 & 57.50 & 24.60 & 18.00 \\
			& 2 & 6744.00 & 8607.80 & 2596.00 & 38.80 & 24.50 & 41.70 \\
			& 3 & 6744.00 & 8607.80 & 2174.00 & 22.10 & 32.00 & 46.00 \\
			& 4 & 6744.00 & 8607.80 & 1489.00 & 18.00 & 24.50 & 47.60 \\
			& 5 & 6744.00 & 8607.80 & 2667.00 & 42.10 & 28.10 & 36.50 \\
			& 6 & 6744.00 & 8607.80 & 2162.00 & 40.38 & 24.72 & 39.98 \\
			& 7 & 6744.00 & 8607.80 & 2545.00 & 40.38 & 27.74 & 40.05 \\
			\midrule
			\text{2} & 1 & 37,635.00 & 17,532.00 & 23,352.00 & 73.00 & 64.00 & 51.90 \\
			& 2 & 37,635.00 & 17,532.00 & 14,668.00 & 78.30 & 72.70 & 76.40 \\
			& 3 & 37,635.00 & 17,532.00 & 3396.00 & 87.70 & 70.50 & 74.40 \\
			& 4 & 37,635.00 & 17,532.00 & 15,453.00 & 73.00 & 74.00 & 64.00 \\
			& 5 & 37,635.00 & 17,532.00 & 13,779.00 & 69.00 & 56.90 & 76.70 \\
			& 6 & 37,635.00 & 17,532.00 & 3694.00 & 73.17 & 85.90 & 74.03 \\
			& 7 & 37,635.00 & 17,532.00 & 2025.00 & 73.17 & 62.49 & 78.14 \\
			
			\bottomrule
\end{tabularx}
   \end{adjustwidth}
    \end{table}
     \begin{table}[H]\ContinuedFloat
\caption{{\em Cont.}} 
    \begin{adjustwidth}{-\extralength}{0cm}
    \newcolumntype{C}{>{\centering\arraybackslash}X}
		\begin{tabularx}{\fulllength}{ccccCccc}
			\toprule
			\textbf{Part Type} & \textbf{Build Orientation} & \textbf{Volume [mm$^3$]} & \textbf{Surface Area [mm$^2]$} & \textbf{Supporting Structure Volume~[mm$^3$]} & \textbf{Length [mm]} & \textbf{Width [mm]} & \textbf{Heigth  [mm]} \\
			\midrule
			\text{3} & 1 & 1029.00 & 1017.00 & 98.00 & 28.30 & 13.80 & 13.80 \\
			& 2 & 1029.00 & 1017.00 & 142.00 & 24.90 & 13.80 & 25.10 \\
			& 3 & 1029.00 & 1017.00 & 376.00 & 21.80 & 15.90 & 26.90 \\
			& 4 & 1029.00 & 1017.00 & 0.00 & 13.70 & 13.80 & 28.30 \\
			& 5 & 1029.00 & 1017.00 & 536.00 & 26.70 & 14.90 & 22.40 \\
			& 6 & 1029.00 & 1017.00 & 155.00 & 24.87 & 13.75 & 24.87 \\
			& 7 & 1029.00 & 1017.00 & 138.55 & 19.77 & 21.95 & 25.14 \\
			\midrule
			\text{4} & 1 & 105,909.00 & 45,458.30 & 36,239.00 & 69.00 & 169.00 & 36.60 \\
			& 2 & 105,909.00 & 45,458.30 & 10,143.00 & 75.70 & 165.70 & 101.00 \\
			& 3 & 105,909.00 & 45,458.30 & 34,743.00 & 139.00 & 93.20 & 155.70 \\
			& 4 & 105,909.00 & 45,458.30 & 44,612.00 & 69.80 & 36.60 & 169.00 \\
			& 5 & 105,909.00 & 45,458.30 & 47,096.00 & 90.60 & 147.40 & 136.20 \\
			& 6 & 105,909.00 & 45,458.30 & 13,352.00 & 85.35 & 158.36 & 122.65 \\
			& 7 & 105,909.00 & 45,458.30 & 20,212.00 & 141.10 & 90.55 & 157.22 \\
			\midrule
			\text{5} & 1 & 28,588.10 & 20,398.80 & 1183.00 & 77.00 & 77.00 & 60.90 \\
			& 2 & 28,588.10 & 20,398.80 & 5542.00 & 77.00 & 77.00 & 60.90 \\
			& 3 & 28,588.10 & 20,398.80 & 25,217.00 & 76.70 & 60.90 & 76.70 \\
			& 4 & 28,588.10 & 20,398.80 & 25,150.00 & 69.70 & 74.60 & 76.70 \\
			& 5 & 28,588.10 & 20,398.80 & 12,042.00 & 73.20 & 74.10 & 71.10 \\
			& 6 & 28,589.10 & 20,398.80 & 25,145.00 & 72.45 & 73.73 & 76.90 \\
			& 7 & 28,589.10 & 20,398.80 & 25,145.00 & 74.02 & 72.47 & 76.74 \\
			\midrule
			\text{6} & 1 & 6310.00 & 9792.00 & 3908.00 & 16.60 & 79.70 & 11.50 \\
			& 2 & 6310.00 & 9792.00 & 3726.00 & 34.70 & 71.20 & 60.60 \\
			& 3 & 6310.00 & 9792.00 & 5187.00 & 27.20 & 37.20 & 80.40 \\
			& 4 & 6310.00 & 9792.00 & 2425.00 & 16.60 & 11.50 & 79.70 \\
			& 5 & 6310.00 & 9792.00 & 3737.00 & 29.60 & 73.20 & 57.40 \\
			& 6 & 6310.00 & 9792.00 & 2439.00 & 30.78 & 16.17 & 81.23 \\
			& 7 & 6310.00 & 9792.00 & 3267.00 & 26.36 & 62.57 & 65.43 \\
			\bottomrule
		\end{tabularx}
	\end{adjustwidth}
\end{table}

\begin{adjustwidth}{-\extralength}{0cm}
\reftitle{References}

\PublishersNote{}
\end{adjustwidth}
\end{document}